\input amstex\documentstyle{amsppt}  
\pagewidth{12.5cm}\pageheight{19cm}\magnification\magstep1
\topmatter
\title Restriction of a character sheaf to conjugacy classes\endtitle
\author G. Lusztig\endauthor
\address{Department of Mathematics, M.I.T., Cambridge, MA 02139}\endaddress
\thanks{Supported in part by National Science Foundation grant DMS-0758262}\endthanks
\endtopmatter   
\document
\define\Irr{\text{\rm Irr}}

\define\si{\sim}

\define\sqc{\sqcup}

\define\hE{\hat E}

\define\hG{\hat G}

\define\hY{\hat Y}

\define\bA{\bar A}

\define\bY{\bar Y}

\define\lb{\linebreak}

\define\op{\oplus}

\define\part{\partial}
\define\em{\emptyset}

\define\n{\notin}

\define\m{\mapsto}
\define\do{\dots}

\define\bsl{\backslash}

\define\lra{\leftrightarrow}

\define\sub{\subset}    
\define\bxt{\boxtimes}
\define\T{\times}
\define\ti{\tilde}
\define\nl{\newline}
\redefine\i{^{-1}}

\define\ov{\overline}
\define\ot{\otimes}
\define\bbq{\bar{\QQ}_l}

\define\End{\text{\rm End}}

\define\ind{\text{\rm ind}}

\define\a{\alpha}
\redefine\b{\beta}

\define\g{\gamma}
\redefine\d{\delta}
\define\e{\epsilon}
\define\et{\eta}

\redefine\o{\omega}
\define\p{\pi}
\define\ph{\phi}

\define\r{\rho}
\define\s{\sigma}
\redefine\t{\tau}

\define\k{\kappa}
\redefine\l{\lambda}
\define\z{\zeta}

\redefine\G{\Gamma}
\redefine\D{\Delta}

\define\Si{\Sigma}

\define\Ph{\Phi}

\define\kk{\bold k}

\define\FF{\bold F}

\define\NN{\bold N}

\define\QQ{\bold Q}

\define\SS{\bold S}

\define\WW{\bold W}
\define\ZZ{\bold Z}

\define\cb{\Cal B}

\define\ce{\Cal E}
\define\cf{\Cal F}

\define\ch{\Cal H}

\define\ck{\Cal K}
\define\cl{\Cal L}

\define\cs{\Cal S}

\define\cu{\Cal U}
\define\cv{\Cal V}
\define\cw{\Cal W}
\define\cz{\Cal Z}

\define\fc{\frak c}

\define\fF{\frak F}

\define\fO{\frak O}

\define\fS{\frak S}

\define\fU{\frak U}

\define\tY{\ti Y}

\define\bP{\bar P}

\define\tce{\ti\ce}

\define\tcw{\ti{\cw}}
\define\tfc{\ti{\fc}}
\define\BM{BM}
\define\DL{DL}
\define\EF{E}
\define\COX{L1}
\define\REP{L2}
\define\ORA{L3}
\define\ICC{L4}
\define\CSI{L5}
\define\CSII{L6}
\define\CSIII{L7}
\define\INT{L8}
\define\SUP{L9}
\define\CDGIII{L10}
\define\CLE{L11}
\define\FAM{L12}

\head Introduction\endhead
\subhead 0.1\endsubhead
Let $\kk$ be an algebraically closed field of characteristic $p\ge0$. Let $G$ be a connected reductive 
group over $\kk$. Let $\hG$ be the set of (isomorphism classes of) character sheaves on $G$. (Recall that 
the character sheaves of $G$ are certain simple perverse sheaves on $G$ (see \cite{\INT}) 
which are equivariant under $G$-conjugation.) In this paper we are interested in studying the restriction of 
a character sheaf of $G$ to a conjugacy class of $G$. For $g\in G$ let $g=g_sg_u=g_ug_s$ be the Jordan 
decomposition of $g$ ($g_s$ is semisimple, $g_u$ is unipotent). For a subset $R$ of $G$ let 
$R_s=\{g_s;g\in X\}$, $R_u=\{g_u;g\in X\}$. The set $\hG$ can be naturally
partitioned into equivalence classes called {\it families} as in \cite{\SUP, 10.6} and (assuming that $p$ 
is not a bad prime for $G$) to each family $\fF$ one can attach a unipotent class $C=C_{\fF}$ of $G$ as in 
\cite{\SUP, 10.5} so that the following property holds (see \cite{\SUP, 10.7}).

(a) {\it If $D$ is a conjugacy class of $G$ such that $\dim D_u\ge\dim C$ and $D_u\ne C$ then $A|_D=0$ for 
any $A\in\fF$. There exists a conjugacy class $D$ of $G$ and $A\in\fF$ such that $D_u=C$ and $A|_D\ne0$.}
\nl
Note that (a) characterizes $C$ uniquely. (Actually, in \cite{\SUP}, $p$ is assumed to be sufficiently large
or $0$ but from this the case where $p$ is only assumed to be not a bad prime can be deduced by standard 
methods.) If $A\in\hG$ belongs to a family $\fF$, the unipotent class $C_{\fF}$ is said to be the 
{\it unipotent support} of $A$.

We now state the following refinement of (a).

\proclaim{Theorem 0.2} Assume that $p$ is not a bad prime for $G$. 
Let $\fF,C=C_{\fF}$ be as above. Let $D$ 
be a conjugacy class of $G$ such that $D_u=C$. Then for any $A\in\fF$ we have $A|_D=\cl[\dim(D)+c]$ where 
$\cl$ is a local system on $D$ and $c\in\NN$ depends only on $A$, not on $D$.
\endproclaim
The proof is given in \S1.

\subhead 0.3\endsubhead
Recall that in \cite{\ICC, 3.1} we have defined a partition of $G$ into finitely many 
locally closed smooth irreducible subvarieties of $G$ invariant by conjugation (called strata).
From the definitions we see that if $Y$ is a stratum of $G$ then $Y_u$ is a single unipotent conjugacy 
class of $G$. We have the following result.

\proclaim{Corollary 0.4}  Assume that $p$ is not a bad prime for $G$. Let $A$ be a character sheaf on $G$ 
and let $Y$ be a stratum of $G$. Let $\fF$ be the family in $\hG$ that contains $A$ and let $C=C_{\fF}$.

(a) If $\dim Y_u\ge\dim C$ and $Y_u\ne C$ then $A|_Y=0$.

(b) If $Y_u=C$ then $A|_Y$ is a local system (up to shift).
\endproclaim
(a) follows immediately from the definition of $C$. We prove (b). Since all conjugacy classes contained in 
$Y$ have the same dimension, we see from Theorem 0.2 that there exists $i\in\ZZ$ such that
$\ch^jA|_Y=0$ for $j\ne i$. It remains to show that $\ch^iA|_Y$ is a local system. This follows from
\cite{\CSIII, 14.2(a)}.

\subhead 0.5\endsubhead
{\it Notation.} Let $\cb$ be the variety of Borel subgroups of $G$. We write 
$\cb\T\cb=\sqc_{w\in\WW}\fO_w$ where $\fO_w (w\in\WW)$ are the orbits of $G$ acting on $\cb\T\cb$ by 
simultaneous conjugation and $\WW$ is the Weyl group (it is naturally
a finite Coxeter group with set of simple reflections $\SS$). 
Let $B^*$ be a fixed Borel subgroup of $G$ and let $T^*$ be a fixed maximal torus of $B^*$. 
The parabolic subgroups of $G$ containing $B^*$ are in natural bijection $J\lra P_J$
with the subsets $J\sub\SS$; thus $P_\em=B^*,P_\SS=G$. For $J\sub\SS$ let $L_J$ be the unique Levi subgroup 
of $P_J$ such that $T^*\sub L_J$.

If $P$ is a parabolic subgroup of $G$ we denote by $U_P$ the unipotent radical of $P$ and by 
$\p_P:P@>>>\bP:=P/U_P$ the canonical homomorphism. For an affine algebraic group $H$ let $H^0$ be the 
identity component of $H$; let $\cz_H$ be the centre of $H$. We denote by $l$ a prime number invertible in
$\kk$. All sheaves (in particular local systems) are assumed to be $\bbq$-sheaves. For a group $\D$ and 
$g\in\D$ let $Z_\D(g)$ be the centralizer of $g$ in $\D$. If $\D$ is finite we denote by $\Irr\D$ a set of 
representatives for the isomorphism classes of irreducible representations of $\D$ over $\bbq$. 

If $\kk$ is an algebraic closure of the finite field $\FF_p$ with with $p$ elements and $q$ is a power of $p$
we denote by $\FF_q$ the subfield of $\kk$ with $q$ elements. 

\subhead 0.6\endsubhead
In \S2 we describe a parametrization of the unipotent character sheaves of $G$ (suggested by Theorem 0.2) in
terms of restrictions to various conjugacy classes assuming, that $p$ is not a bad prime. This 
description fails in bad characteristic.

In \S3 we establish a canonical bijection between the set of unipotent character sheaves on $G$ and the
set of unipotent representations of a split reductive group of the same type over a finite field (here there
is no restriction on $p$). In fact both sets are put in canonical bijection with a combinatorially defined
set $\fS_{\WW}$ defined purely in terms of the Weyl group $\WW$.

\subhead 0.7 \endsubhead
Now assume that $G$ has connected centre and let $\cl$ be a local system of rank $1$ on $T^*$ such that
$\cl^{\ot m}\cong \bbq$ for some $m\ge1$ invertible in $\kk$. Let $\hG_{\cl}$ be the subset of $\hG$
defined as in \cite{\INT, 4.2}. Then the method in \S3 can be extended to give a canonical bijection of 
$\hG_{\cl}$ with $\fS_{\WW_{\cl}}$ where $\WW_{\cl}$ is the stabilizer of $\cl$ in $\WW$.

\head 1. Proof of Theorem 0.2\endhead
\subhead 1.1\endsubhead
A cuspidal pair for $G$ is a pair $(\Si,\ce)$ where $\Si$ is the inverse image under $G@>>>G/\cz_G^0$ of an 
isolated conjugacy class in $G/\cz_G^0$ and $\ce$ is a local system on $\Si$ such that for some $n\ge1$ 
invertible in $\kk$,

$\ce$ is equivariant for the $G\T\cz_G^0$ action $(g,z):g_1\m z^ngg_1g\i$ on $\Si$;

for any parabolic subgroup $P\ne G$ and any $y\in\bP$, we have $H^\d_c(\p_P\i(y)\cap\Si,\ce)=0$ where 
$\d=\dim(\Si/\cz_G^0)-\dim(\bP)+\dim Z_{\bP}(y)$.
\nl
(See \cite{\ICC, \S2}).

An induction datum for $G$ is a triple $(L,\Si,\ce)$ where $L$ is a Levi subgroup of a parabolic subgroup 
of $G$ and $(\Si,\ce)$ is a cuspidal pair for $L$. Given an induction datum $(L,\Si,\ce)$ for $G$ we set

$\Si_r=\{g\in\Si;Z_G(g_s)^0\sub L\}$

$Y=\sqc_{x\in G}x\Si_rx\i$

$\tY=\{(g,xL)\in G\T G/L;x\i gx\in\Si_r\}$

$\hY=\{(g,x)\in G\T G;x\i gx\in\Si_r\}$.
\nl
We have a diagram
$$\Si@<\a<<\hY@>\b>>\tY@>\p>>Y$$
where $\a(g,x)=x\i gx$, $\b(g,x)=(g,xL)$, $\p(g,xL)=g$. The local system $\a^*\ce$ on $\hY$ is 
$L$-equivariant for the action of $L$ on $\hY$ given by $l:(g,x)\m(g,xl\i)$ hence it is equal to $\b^*\tce$ 
for a well defined 
local system $\tce$ on $\tY$. Now $\p_!\tce$ is a well defined local system on $Y$ (a locally closed smooth 
irreducible subvariety of $G$). Let $K=IC(\bY,\p_!\tce)$ extended to $G$ by $0$ on $G-\bY$.
Here $\bY$ is the closure of $Y$ in $G$. We set $f_0=\dim Y=\dim G-\dim L+\dim\Si$. Note that $K[f_0]$ is a 
perverse sheaf.

\subhead 1.2\endsubhead
Let $\s\in G$ be a semisimple element which is $G$-conjugate to an element of $\Si_s$. Let 
$M=\{x\in G;x\i\s x\in\Si_s\}$. We have $M\ne\em$. Let $\G=Z_G^0(\s)\bsl M/L$, a finite nonempty set. The
group $\tcw_\Si:=\{y\in G;yLy\i=L,y\Si y\i=\Si\}$ acts on $\G$ by $y:\et\m\et y\i$. This induces an action 
of the finite group $\cw_\Si=\tcw_\Si/L$ on $\G$.

For any $x\in M$ let $L_x=xLx\i\cap Z^0_G(\s)$; this is a Levi subgroup of a parabolic subgroup of 
$Z_G^0(\s)$. Let $C_x=\{v\in Z^0_G(\s);v\text{ unipotent },x\i\s vx\in\Si\}$; this is a unipotent class in 
$Z_G^0(\s)$, see \cite{\CSII, 7.11(c)}. Let $\Si_x=\cz_{L_x}^0C_x$. Let $\ce_x$ be the local system on 
$\Si_x$ obtained as the inverse image of $\ce$ under the map $\Si_x@>>>\Si$, $g\m x\i\s gx$. Then 
$(\Si_x,\ce_x)$ is a cuspidal pair for $L_x$. (See \cite{\CSII, 7.11(a)}.) The definitions in 1.1 are 
applicable to $Z_G^0(\s),L_x,\Si_x,\ce_x$ instead of $G,L,\Si,\ce$. Let 
$\p_x:tY'_x@>>>Y'_x,\bY'_x,\tce_x,K_x$ be obtained from $\p:\tY@>>>Y,\bY,\tce,K$ in 1.1 by replacing 
$G,L,\Si,\ce$ by $Z_G^0(\s),L_x,\Si_x,\ce_x$. We set $f=\dim Y'_x=\dim Z_G^0(\s)-\dim L_x+\dim\Si_x$; this
is independent of the choice of $x$ in $M$ (it is in fact equal to $\dim Z_G^0(\s)-\dim L+\dim\Si$, see the 
proof of \cite{\CDGIII, 16.10(b)}). Note that $K_x[f]$ is a perverse sheaf on $Z_G^0(\s)$ with 
support equal to $\bY'_x$.

For any $\et\in\G$ (viewed as a subset of $M$) we choose a base point $x_\et\in\et$. We set 
$$\align&L_\et=L_{x_\et},C_\et=C_{x_\et},\Si_\et=\Si_{x_\et},\p_\et=\p_{x_\et},\tY'_\et=\tY'_{x_\et},\\&
Y'_\et=Y'_{x_\et},\bY'_\et=\bY'_{x_\et},\ce_\et=\ce_{x_\et},\tce_\et=\tce_{x_\et},K_\et=K_{x_\et}.\endalign$$

\subhead 1.3\endsubhead
Let $\fU$ be the set of all open sets $\cu$ in $Z^0_G(\s)$ such that

$1\in\cu$;

$g\cu g\i=\cu$ for all $g\in Z^0_G(\s)$;

for $x\in Z_G^0(\s)$ we have $x\in\cu$ if and only if $x_s\in\cu$.
\nl
For example, $Z^0_G(\s)\in\fU$.

Let $\cu\in\fU$. For $\et\in\G$ we set 

$\tY_{\cu,\et}=\{(g,xL)\in\tY;g\in\s\cu,x\in\et\}$, $Y_{\cu,\et}=\p(\tY_{\cu,\et})$.
\nl
Let $y\in\tcw_\Si,\et\in\G$. If $(g,xL)\in\tY_{\cu,\et}$ then $(g,xy\i L)\in\tY_{\cu,\et y\i}$ hence
$Y_{\cu,\et y\i}=Y_{\cu,\et}$. Thus for any $\cw_\Si$-orbit $Z$ in $\G$ we can set $Y_{\cu,Z}=Y_{\cu,\et}$
where $\et$ is any element in $Z$. For $Z$ as above we set 
$\cv_Z=\cap_{\et\in Z}(\s(Y'_\et\cap\cu)\cap Y_{\cu,Z})$.
Let $\cv=\cup_Z\cv_Z$ where $Z$ runs over the $\cw_\Si$-orbits in $\G$. For $\et\in\G$ let $[\et]$ be the
$\cw_\Si$-orbit of $\et$ in $\G$. 

As shown in \cite{\CDGIII,\S16} (see also \cite{\CSII, \S8}), for any $\cu_1\in\fU$ we can find 
$\cu\in\fU$ with $\cu\sub\cu_1$ so that (a),(b),(c) below hold.

(a) $\cv$ is an open subset of $Y\cap\s\cu$ of pure dimension $f$ and the $\cv_Z$ ($Z$ runs over the 
$\cw_\Si$-orbits in $\G$) form a finite partition of $\cv$ into open and closed subsets; $\cv$ is open 
dense in $\bY\cap\s\cu=\cup_{\et\in\G}\s(\bY'_\et\cap\cu)$; for $\et\in\G$, $\s\i\cv_{[\et]}$ is open dense 
in $\bY'_\et\cap\cu$.

(b) Let ${}^0\tY=\{(g,xL)\in\tY;g\in\cv\}$. For $\et\in\G$ let 
${}^0\tY'_\et=\{(h,zL_\et)\in\tY'_\et;h\in\s\i\cv_{[\et]}\}$. We have a commutative diagram
$$\CD\sqc_{\et\in\G}{}^0\tY'_\et@>a>>{}^0\tY\\
@VVV                         @VVV\\
\s\i\cv @>\e>>                \cv\endCD$$
where $a$ is the isomorphism given by $(h,zL_\et)\m(\s g,zx_\et L)$, $\e(h)=\s h$ and the vertical maps are 
given by the first projection.

(c) The canonical isomorphism
$$\op_{\et\in\G}(\p_{\et!}\tce_\et)|_{\s\i\cv_{[\et]}}@>\si>>\e^*((\p_!\tce)|_\cv)$$
of local systems over $\s\i\cv$ obtained from (b) (where $(\p_{\et!}\tce_\et)|_{\s\i\cv_{[\et]}}$ 
is extended by $0$ on $\s\i\cv_{Z'}$ for $Z'\ne[\et]$) extends uniquely to an isomorphism 
$$\k:\op_{\et\in\G}K_\et|_\cu@>\si>>\e'{}^*(K|_{\s\cu})$$
where $\e':\cu@>>>\s\cu$ is $g\m\s g$.
\nl
(The uniqueness follows from the definition of the intersection cohomology complex.)

\proclaim{Proposition 1.4} Let $E$ be a semisimple class of $G$. Let $\hE=\{g\in G;g_s\in E\}$. Let 
$L,\Si,\ce,Y,\bY,K,f_0$ be as in 1.1. Let $c=\dim\cz_L^0$. Then $K|_{\hE}[f_0-c]$ is a semisimple 
perverse sheaf.
\endproclaim
Let $\s\in E_s$. If $\s\n\bY_s$ then clearly $K|_{\hE}=0$. Thus we can assume that $\s\in\bY_s$. It follows 
that $\s$ is as in 1.2. Let $\Xi$ be the set of unipotent elements in $Z_G^0(\s)$. Note that $\s\Xi\sub\hE$.
It is enough to show that $K[f-c]|_{\s\Xi}$ is a semisimple perverse sheaf where $f$ is as in 1.2. (Note 
that $\dim(G/Z_G^0(\s))=f_0-f$.) With notation of 1.3(c) we have 
$$\op_{\et\in\G}K_\et|_\Xi@>\si>>\e''{}^*(K|_{\s\Xi})\tag a$$
where $\e'':\Xi@>>>\s\Xi$ is $g\m\s g$. (We use that $\Xi\sub\cu$.) We see that it is enough to show that 
for any $\et\in\G$, $K_\et|_\Xi[f-c]$ is a semisimple perverse sheaf.
We define a local system $\ce_\et^1$ on $\Si_\et=\cz_{L_\et}^0C_\et=\cz_L^0C_\et$ ($C_\et$ as in 1.2) as 
$\bbq\bxt(\ce_\et|_{1\T C_\et})$ where $\bbq$ is viewed as a local system on $\cz_L^0$. We define $K^1_\et$ 
in terms of $Z_G^0(\s),L_\et,\Si_\et,\ce_\et^1$ in the same way as $K$ was defined in terms of $G,L,\Si,\ce$
in 1.1.From \cite{\CDGIII, 15.2} we see that there is a canonical isomorphism
$$(K_\et)|_\Xi@>\si>>(K^1_\et)|_\Xi.\tag b$$
It is then enough to show that $K^1_\et|_\Xi[f-c]$ is a semisimple perverse sheaf. But this follows from 
\cite{\ICC, (6.6.1)}. The proposition is proved.

\subhead 1.5\endsubhead
{\it Remark.} In the special case where $E=\{1\}$, $L=\Si$ is a maximal torus of $G$, $\ce=\bbq$, the 
proposition above reduces to a result in \cite{\BM}.

\proclaim{Corollary 1.6} Let $A$ be a character sheaf of $G$ and let $E,\hE$ be as in 1.4. Then for some 
integer $m$, $A|_{\hE}[m]$ is a semisimple perverse sheaf.
\endproclaim
We can find $L,\Si,\ce,Y,\bY,K$ as in 1.1 such that $A$ is a direct summand of $K[f_0]$ ($f_0$ as in 1.1). 
Hence $A|_{\hE}[-c]$ is a direct summand of $K|_{\hE}[f_0-c]$ which by 1.4 is a semisimple perverse sheaf. 
It follows that $A|_{\hE}[-c]$ is itself a semisimple perverse sheaf. 

\subhead 1.7\endsubhead
We prove Theorem 0.2. Let $A\in\fF$. Let $C=C_{\fF}$. We can find $L,\Si,\ce,Y,\bY,K$ as in 1.1 such that 
$A$ is a direct summand of $K[f_0]$ ($f_0$ as in 1.1). We can also assume that $\ce$ is irreducible.
Let $\s\in D_s$. If $\s\n\bY_s$ then clearly $K|_D=0$ hence $A|_D=0$. Thus we can assume that $\s\in\bY_s$. 
It follows that $\s$ is as in 1.2. Let $\Xi$ be the set of unipotent elements in $Z_G^0(\s)$. 
By the proof of Proposition 1.4 and with the notation there we have
$$\e''{}^*(K|_{\s\Xi})[f-c]\cong\op_{\et\in\G}K^1_\et[f-c]|_\Xi\cong A_1\op A_2\op\do\op A_n$$
where $A_1,\do,A_n$ are simple perverse sheaves on $\Xi$, equivariant under conjugation by $Z_G^0(\s)$.
Since $\e''{}^*(A|_{\s\Xi})[f-c-f_0]$ is a direct summand of $\e''{}^*(K_{\s\Xi})[f-c]$ it follows that
$\e''{}^*(A_{\s\Xi})[f-c-f_0]$ is isomorphic to a direct summand of $A_1\op A_2\op\do\op A_n$ hence
$$\e''{}^*(A|_{\s\Xi})[f-c-f_0]\cong\op_{j\in J}A_j\text{ for some }J\sub[1,n].\tag a$$
Now for each $j\in J$ there is a unique unipotent class $\fc_j$ of $Z_G^0(\s)$
such that $A_j|_{\fc_j}$ is of the form $\cl_j[\dim\fc_j]$ where $\cl_j$ is an irreducible local
system and $A_j|_{\fc'}=0$ for any unipotent class $\fc'$ of $Z_G^0(\s)$ such that $\fc'\not\sub\ov{\fc_j}$.
It follows that for any unipotent class $\fc$ of $Z_G^0(\s)$ we have
$$\e''{}^*(A|_{\s\fc})[f-c-f_0]\cong\op_{j\in J;\fc\sub\ov{\fc_j}-\fc_j}A_j|_\fc\op\op_{j\in J;\fc=\fc_j}
\cl_j[\dim\fc].\tag b$$
Now assume that $\fc$ is a unipotent class of $Z_G^0(\s)$ such that $\fc\sub C$. Assume that 

(c) $\op_{j\in J; \fc\sub\ov{\fc_j}-\fc_j}A_j|_\fc\ne0$
\nl
that is, there exists $j_0\in J$ such that $\fc\sub\ov{\fc_{j_0}}-\fc_{j_0}$, $A_{j_0}|_\fc\ne0$. Using (a) 
we deduce
$$\e''{}^*(A|_{\s\fc_{j_0}})[f-c-f_0]\cong\op_{j\in J}A_j|_{\fc_{j_0}}.$$
The last direct sum contains as a summand $A_{j_0}|_{\fc_{j_0}}=\cl_{j_0}[\dim\fc_{j_0}]\ne0$ hence the whole
direct sum is $\ne0$ and $\e''{}^*(A|_{\s\fc_{j_0}})[f-c-f_0]\ne0$. Thus $A|_{\s\fc_{j_0}}\ne0$. Let $\tfc$
be the unipotent class of $G$ that contains $\fc_{j_0}$. Since $\fc\sub\ov{\fc_{j_0}}-\fc_{j_0}$, we see 
that $C\sub\ov{\tfc}-\tfc$.
If $D'$ is the conjugacy class of $G$ that contains $\s\fc_{j_0}$ we have $D'_u=\tfc$, $\dim D'_u>\dim C$.
Hence by by 0.1(a) we have $A_{D'}=0$ so that $A|_{\s\fc_{j_0}}=0$, a contradiction. Thus the assumption (c) 
cannot hold so that the direct sum in (c) is zero and (b) reduces to
$$\e''{}^*(A|_{\s\fc})[f-c-f_0]\cong\op_{j\in J;\fc=\fc_j}\cl_j[\dim\fc].$$
We see that $A|_{\s\fc}[f-c-f_0-\dim\fc]$ is a local system.
Applying this to $\fc$ equal to one of $\fc_1,\fc_2,\do,\fc_r$, the unipotent classes in $Z_G^0(\s)$ that 
are contained in $C$ and satisfy $\{g\in D;g_s=\s\}=\s\fc_1\sqc\s\fc_2\sqc\do\sqc\s\fc_r$ we deduce 
that $A|_{\s\fc_1\sqc\do\sqc\s\fc_r}[f-c-f_0-d]$ is a local system.
(Here $d=\dim\fc_1=\do=\dim\fc_r$; note that $\fc_1,\do\fc_r$ are permuted transitively by $Z_G(\s)$ hence 
have the same dimension). It also follows that $A|_D[f-c-f_0-d]$ is a local system. 
We have $f_0-f+d=\dim(G/Z_G^0(\s))+d=\dim(D)$. Hence $A|_D[-\dim(D)-c]$ is a local system. Theorem 0.2 is 
proved.

\subhead 1.8\endsubhead
The following result can be proved in the same way as Theorem 0.2 with no restriction on $p$.

{\it Let $A\in\hG$. Assume that $C_0$ is a unipotent class of $G$ such that the following holds: if $D$ is a
conjugacy class of $G$ such that $\dim D_u\ge\dim C_0$ and $D_u\ne C_0$ then $A|_D=0$. Then for any 
conjugacy class $D$ of $G$ such that $D_u=C$ we have $A|_D=\cl[\dim(D)+c]$ where $\cl$ is a local system on 
$D$ and $c\in\NN$ depends only on $A$, not on $D$.}

\subhead 1.9\endsubhead
The local system $\cl$ which appears in Theorem 0.2 can be reducible. For example if $G$ is of type $B_2$ 
(resp. $G_2$) and $p$ is not a bad prime for $G$ then there exist character sheaves on $G$ with unipotent 
support equal to the subregular unipotent class $C$ in $G$ whose restriction to $C$ is (up to shift) the 
direct sum of $\bbq$ and an irreducible local system of rank $1$ (resp. $2$).

\head 2. Unipotent character sheaves\endhead
\subhead 2.1\endsubhead
For $w\in\WW$ let $Y_w=\{(g,B)\in G\T\cb;(B,gBg\i)\in\fO_w\}$ and let $\p_w:Y_w@>>>G$ be the first 
projection. Let $\hG^{un}$ be the subset of $\hG$ consisting of unipotent character sheaves that is, those
character sheaves on $G$ which appear as constituents of the perverse cohomology sheaf 
${}^pH^i(\p_{w!}\bbq)$ for some $w\in\WW,i\in\ZZ$. Note that $\hG^{un}$ is a union of families of $\hG$. In 
fact, as in \cite{\INT, 4.6}, we have a partition $\hG^{un}=\sqc_{\cf}\hG^{un}_{\cf}$ where $\cf$ runs over 
the families in $\Irr\WW$ and for each family $\cf$ of $\Irr\WW$, $\hG^{un}_{\cf}$ is a family of $\hG$.

If $A\in\hG^{un}$ the proof of Theorem 0.2 for $A$ simplifies somewhat. In this case Proposition 1.4 is only 
needed in the case where $(L,C,\ce)$ satisfies the condition that $\ce$ is equivariant for the
$L\T\cz_L^0$-action $(g,z):g_1\m zgg_1g\i$ on $\Si$ (hence $\ce$ is the inverse image under 
$\Si@>>>\Si/\cz_L^0$ of an irreducible local system on $\Si/\cz_L^0$). Hence if $\s$ is as in 1.2 then with 
notation in the proof of 1.4, for $\et\in\G$, the local system $\ce_\et$ on $\Si_\et$ is equivariant for the
action of $\cz_{L_\et}^0=\cz_L^0$ on $\Si_\et$ (left multiplication) hence $\ce_\et=\ce_\et^1$ and 
$K^1_\et=K_\et$. Thus the step 1.4(b) in the proof of 1.4 is unnecessary in this case.

We shall try to make Theorem 0.2 more precise in the case of unipotent character sheaves.
In the remainder of this section we assume that $p$ is not a bad prime for $G$.

\subhead 2.2\endsubhead
For a finite group $\D$, $M(\D)$ is the set of all pairs $(g,\r)$ where $g\in\D$ is defined up to conjugacy 
and $\r\in\Irr Z_\D(g)$.

\subhead 2.3\endsubhead
In the remainder of this section we fix a family $\cf$ of $\Irr\WW$. Let $\fF=\hG^{un}_{\cf}$ be the 
corresponding family in $\hG^{un}$. In this case $C=C_{\fF}$ is the special unipotent class of $G$ such that
the corresponding Springer representation of $\WW$ belongs to $\cf$.
Let $\cs_C$ be the set of conjugacy classes $D$ in $G$ such that $D_u=C$.
Let $\o\in C$. Let $A(\o)=Z_G(\o)/Z_G^0(\o)$ and let $\D=\bA(\o)$ be the canonical quotient of $A(\o)$
defined in \cite{\ORA, (13.1.1)}. Let $Z_G(\o)@>j'>>A(\o)@>h>>\D$ be the obvious 
(surjective) homomorphisms; let $j=hj':Z_G(\o)@>>>\D$. Let $[\D]$ be the set of conjugacy classes in $\D$. 
For $D\in\cs_C$ let $\ph(D)$ be the conjugacy class of $j(g_s)$ in $\D$ where $g\in D$ is such that 
$g_u=\o$; 
clearly such $g$ exists and is unique up to $Z_G(\o)$-conjugacy so that the conjugacy class of $j(g_s)$ is 
independent of the choice of $g$. Thus we get a (surjective) map $\ph:\cs_C@>>>[\D]$. For $\g\in[\D]$ we set 
$\cs_{C,\g}=\ph\i(\g)$. We now select for each $\g\in[\D]$ an element $x_\g\in\g$.
Let $D\in\cs_{C,\g}$, $\ce\in\Irr Z_\D(x_\g)$. We can find 
$g\in D$ such that $g_u=\o,j(g_s)=x_\g$ (and another choice for such $g$ must be of the form $bgb\i$ where 
$b\in Z_G(\o)$, $j(b)\in Z_\D(x_\g)$). Let $\ce^D$ be the $G$-equivariant local system on $D$ whose stalk at 
$g_1\in D$ is $\{z\in G;zgz\i=g_1\}\T\ce$ modulo the equivalence relation $(z,e)\si(zh\i,j(h)e)$ for all 
$h\in Z_G(g)$. If $g$ is changed to $g_1=bgb\i$ ($b$ as above) then $\ce^D$ is changed to the 
$G$-equivariant local system $\ce^D_1$ on $D$ whose stalk at $g'\in D$ is $\{z'\in G;z'g_1z'{}\i=g'\}\T\ce$
modulo the equivalence relation $(z',e')\si(z'h'{}\i,j(h')e)$ for all $h'\in Z(g_1)$. We have an 
isomorphism of local systems $\ce^D@>\si>>\ce^D_1$ which for any $g'\in D$ maps the stalk of $\ce^D$ at $g'$
to the stalk of $\ce^D_1$ at $g'$ by the rule $(z,e)\m(zb\i,j(b)e)$. (We have $zb\i g_1bz\i=zgz\i=g'$.) This
is compatible with the equivalence relations. Thus the isomorphism class of the local system $\ce^D$ does 
not depend on the choice of $g$. 

Using the methods sketched in \cite{\FAM, \S4} one can prove the following refinement of Theorem 0.2.

\proclaim{Theorem 2.4} (a) Let $A\in\hG^{un}_{\cf}$. There exists a unique $\g\in[\D]$ and a unique 
$\ce\in\Irr Z_\D(x_\g)$ such that

(i) if $D\in\cs_{C,\g}$, we have $A|_D\cong\ce^D[\dim(D)+c]$ ($c$ as in 1.7).

(ii) if $D\in\cs_{C,\g'}$ with $\g'\in[\D]-\{\g\}$, we have $K|_D=0$.

(b) The assignment $A\m(\g,\ce)$ in (a) defines a bijection $\hG^{un}_{\cf}@>\si>>M(\D)$.
\endproclaim

\head 3. Parametrization of unipotent representations and of unipotent character sheaves\endhead
\subhead 3.1\endsubhead   
Let $(W,S)$ be a Weyl group. Let $n=|S|$. Let $\nu_W$ be the number of reflections of $W$. For any subset 
$J$ of $S$ we denote by $W_J$ the subgroup of $W$ generated by $J$.
Assuming that $W$ is irreducible or $\{1\}$ we define a set $\fS^0_W$ as follows.

If $W=\{1\}$ we have $\fS^0_W=\{1\}$.  

If $W$ is of type $A_n (n\ge1)$ we have $\fS^0_W=\em$.  

If $W$ is of type $B_n$ or $C_n$ $(n\ge2)$ we have $\fS^0_W=\{(-1)^{n/2}\}$ if $n=k^2+k$ for some integer 
$k\ge1$ and $\fS^0_W=\em$, otherwise.

If $W$ is of type $D_n$  $n\ge4$ we have $\fS^0_W=\{(-1)^{n/4}\}$ if $n=4k^2$ for some integer $k\ge1$ and
$\fS^0_W=\em$, otherwise.

If $W$ is of type $E_6$ then $\fS^0_W=\{\z\in\bbq;(\z^3-1)/(\z-1)=0\}$.

If $W$ is of type $E_7$ then $\fS^0_W=\{\z\in\bbq;(\z^4-1)/(\z^2-1)=0\}$.

If $W$ is of type $E_8$ then $\fS^0_W=\{\z\in\bbq;(\z^4-1)(z^5-1)(z^6-1)/(\z^2-1)=0\}$ with $\z=1$ 
appearing twice as $\z=1'$ and $\z=1''$.

If $W$ is of type $F_4$ then $\fS^0_W=\{\z\in\bbq;(\z^2-1)(z^3-1)(z^4-1)/(\z^2-1)=0\}$ with $\z=1$ 
appearing twice as $\z=1'$ and $\z=1''$.

If $W$ is of type $G_2$ then $\fS^0_W=\{\z\in\bbq;(\z^2-1)(\z^3-1)/(\z-1)=0\}$.
\nl
The exponents $4,5,6$ which appear in type $E_8$ satisfy: $4=n/2$,\lb  $5\T6=\text{Coxeter number of }W$, 
$4\T5\T6=\nu_W$. 

The exponents $2,3,4$ which appear in type $F_4$ satisfy: $2=n/2$,\lb 
$3\T4=\text{Coxeter number of }W$, $2\T3\T4=\nu_W$.

If $W$ is not irreducible or $\{1\}$, we define $\fS^0_W=\fS^0_{W_1}\T\do\T\fS^0_{W_r}$ where
$W=W_1\T\do\T W_r$ with $W_1,\do,W_r$ ireducible Weyl groups ($r\ge2$).

Assume now that $W$ is irreducible or $\{1\}$. If $J\sub S$ and $\fS^0_{W_J}\ne\em$ then $W_J$ is 
irreducible or $\{1\}$ and for any $J'\sub S$ 
such that $J\sub J'$, conjugation by the longest element $w^{J'}_0$ of $W_{J'}$ leaves $J$ stable; using 
\cite{\COX, 5.9} it follows that the involutions $\s_h:=w^{J\cup h}_0w^J_0=w^J_0w^{J\cup h}$ ($h\in S-J$) 
generate a subgroup of $W$ which is itself a Weyl group denoted by 
$W^{S/J}$. For $h\ne h'$ in $S-J$ the order of $\s_h\s_{h'}$ is equal to 
$$2(\nu_{W_{J\cup h\cup h'}}-\nu_{W_J})/(\nu_{W_{J\cup h}}+\nu_{W_{J\cup h'}}-2\nu_{W_J}).$$
$W^{S/J}$ is a complement to $W_J$ in the normalizer of $W_J$ in $W$.) Let $\fS_W$ be the set of all triples 
$(J,\e,\z)$ where $J\sub S$, $\e\in\Irr W^{S/J}$, $\z\in\fS^0_{W_J}$. We have an imbedding 
$\fS^0_W@>>>\fS_W$, $\z\m(S,1,\z)$. (Note that $W^{S/S}=\{1\}$.)

If $W$ is not irreducible or $\{1\}$, we define $\fS_W=\fS_{W_1}\T\do\T\fS_{W_r}$ where
$W=W_1\T\do\T W_r$ with $W_1,\do,W_r$ ireducible Weyl groups ($r\ge2$).

\subhead 3.2\endsubhead
Now assume that $\kk$ is an algebraic closure of the finite field $\FF_p$ with $p$ elements and that 
$G$ has a fixed $\FF_q$-split rational structure with Frobenius map $F:G@>>>G$. We also assume that 
$F(B^*)=B^*$, $F(T^*)=T^*$.
We fix a square root $\sqrt q$ of $q$ in $\bbq$. Now $F$ induces
a map $F:\cb@>>>\cb$. For any $w\in\WW$ let $X_w=\{B\in\cb;(B,F(B))\in\fO_w\}$ (see \cite{\DL}). Now $G^F$ 
acts on $X_w$ by conjugation hence there is an induced action of $G^F$ on $H^i_c(X_w,\bbq)$ for $i\in\ZZ$. 
Also $X_w$ is stable under 
$F:\cb@>>>\cb$ hence the linear map $F^*:H^i_c(X_w,\bbq)@>>>H^i_c(X_w,\bbq)$ is well defined for any $i$. 
For any $\mu\in\bbq^*$ let
$H^i_c(X_w,\bbq)_\mu$ be the generalized $\mu$-eigenspace of this linear map; it is a $G^F$-submodule of 
$H^i_c(X_w,\bbq)$.
Let $\fU_q$ be the set of all $\r\in\Irr G^F$ such that $\r$ appears in the $G^F$-module $H^i_c(X_w,\bbq)$ 
for some $w\in\WW,i\in\ZZ$. (This is the set of unipotent representations of $G^F$.)
For any $\r\in\fU_q$ and any $w\in\WW$ we denote by $(\r:R_w)$ the multiplicity of $\r$ in the virtual
representation $\sum_i(-1)^iH^i_c(X_w,\bbq)$. 

According to \cite{\REP, 3.9} for any $\r\in\fU_q$ there is a well defined coset $\ti\t_\r$
of $\bbq^*$ modulo its subgroup $\{q^r;r\in\ZZ\}$
such that whenever $\r$ appears in the $G^F$-module $H^i_c(X_w,\bbq)_\mu$ 
(with $i\in\ZZ,w\in\WW,\mu\in\bbq^*$) we have $\mu\in\ti\t_\r$; now $\ti\t_\r$ is contained in a unique coset
(denoted by $\t_\r$) of $\bbq^*$ modulo $\{\sqrt{q}^r;r\in\ZZ\}$.

Let $\fU^0_q$ be the set of all $\r\in\fU_q$ which are cuspidal.
We have the following result. (In the case where $\WW$ is of type $E_8$ (resp. $F_4$) we denote by $w_*$ an 
element of $\WW$ whose characteristic
polynomial in the reflection representation $\WW$ is $(X^4-X^2+1)^2$ (resp. $(X^2-X+1)^2$).)

\proclaim{Theorem 3.3} Assume that $G/\cz_G$ is simple or $\{1\}$. There exists a unique bijection 
$\fS_{\WW}^0@>\sim>>\fU_q^0$ with the following properties. If $\r\in\fU_q^0$ corresponds to 
$\z\in\fS_{\WW}^0$ then $\z\in\t_\r$. If in addition $G/\cz_G$ is of 
type $E_8$ or $F_4$ and $\r'$ (resp. $\r''$) in $\fU_q^0$ corresponds to $1'$ (resp. $1''$) in $\fS_{\WW}^0$
then $(\r':R_{w_*})=1$ and $(\r'':R_{w_*})=0$.
\endproclaim
This follows from \cite{\ORA, 11.2} and its proof.

\subhead 3.4\endsubhead
Now let $J\sub\SS$ and let $\r_0$ be a unipotent cuspidal representation of $L_J^F$. Then 
$L_J/\cz_{L_J}$ is simple or $\{1\}$ hence by 3.3 applied to $L_J$, $\r_0$ corresponds to an element
$\z\in\fS_{\WW_J}^0$. Let $I(J,\z)$ be the representation of $G^F$ induced by $\r_0$ viewed as a 
representation of $P_J^F$.
This is a direct sum of irreducible representations in $\fU_q$; the set of all $\r\in\fU_q$ which appear 
in $I(J,\z)$ is denoted by $\fU_{q,J,\z}$. According to \cite{\REP, 3.26}, the set $\fU_{q,J,\z}$ is in 
natural bijection with the set of irreducible representations of 
a Hecke algebra (a deformation of the Weyl group $\WW^{\SS/J}$ in 3.1 with parameters being powers of $q$ 
explicitly described in \cite{\REP, p.35}); hence it is in natural bijection with the set $\Irr\WW^{\SS/J}$ 
(here we use our choice of $\sqrt q$). We have the following result.

\proclaim{Corollary 3.5} Assume that $G/\cz_G$ is simple or $\{1\}$. There exists a unique bijection 
$\fS_{\WW}@>\si>>\fU_q$ with the following property. If $\r\in\fU_q$ corresponds to $(J,\e,\z)\in\fS_{\WW}$ 
then $\r\in\fU_{q,J,\z}$ and $\r$ corresponds as above to $\e\in\Irr\WW^{\SS/J}$.
\endproclaim

\subhead 3.6\endsubhead
We now drop the assumption on $\kk$ made in 3.2.
Let $A$ be a simple perverse sheaf on $G$ which is equivariant for the conjugation $G$-action on $G$. We 
define an invariant $\l_A\in\bbq^*$ (a root of $1$) as follows. We can find an open dense subset $N$ of the 
support of $A$ 
which is invariant under conjugation by $G$ and an irreducible $G$-equivariant local system $\cl$ on $N$ 
such that $A|_N=\cl[\d]$, $\d=\dim N$. If $g\in N$ then $Z_G(g)$ acts naturally and irreducibly on the stalk 
$\cl_g$. Since $g$ is in the centre of $Z_G(g)$, it acts on $\cl_g$ as a nonzero scalar which is independent
of the choice of $g$ and is denoted by $\l_A$. In particular $\l_A$ is defined for any character sheaf $A$ 
on $G$. (I have found this definition of $\l_A$ in the late 1980's (unpublished); it was also found later in
\cite{\EF}. The definition makes sense also when $G$ is replaced by a finite group $\G$ and $A$ is replaced 
by an irreducible vector bundle on $\G$ equivariant for the conjugation $\G$-action; in this case the 
definition appears in \cite{\ORA, 11.1} and it inspired my later definition for character sheaves.)

Let $\hG^{0,un}$ be the set of all character sheaves in $\hG^{un}$ (see 2.1) which are cuspidal.
For $A\in\hG^{un}$ and $w\in\WW$ we set 
$$(A:\ck_w)=\sum_{i\in\ZZ}(-1)^{\dim G+i}(\text{multiplicity of $A$ in }{}^pH^i(\p_{w!}\bbq))$$
where $\p_{w!}$ is as in 2.1. We have the following result (with $w_*$ as in 3.2).

\proclaim{Theorem 3.7} Assume that $G/\cz_G$ is simple or $\{1\}$. There exists a unique bijection 
$\fS_{\WW}^0@>\sim>>\hG^{0,un}$ with the following properties. If $A\in\hG^{0,un}$ corresponds to 
$\z\in\fS_{\WW}^0$ then $\z=\l_A$. If in addition $G/\cz_G$ is of type $E_8$ or $F_4$ and $A'$ (resp. $A''$)
in $\hG^{0,un}$ corresponds to $1'$ (resp. $1''$) in $\fS_{\WW}^0$ then $(A':\ck_{w_*})=1$ and 
$(A'':\ck_{w_*})=0$.
\endproclaim
This follows from the known classification of unipotent cuspidal character sheaves and the known formulas 
for the numbers $(A:\ck_w)$ (when $p=2$, the cleanness results of \cite{\CLE} must be used).

\subhead 3.8\endsubhead
Now let $J\sub\SS$ and let $A_0$ be a unipotent cuspidal character sheaf on $L_J$. Then 
$L_J/\cz_{L_J}$ is simple or $\{1\}$ hence by 3.7 applied to $L_J$, $A_0$ corresponds to an element
$\z\in\fS_{\WW_J}^0$. The semisimple perverse sheaf $\ind_{P_J}^GA_0$ on $G$ (see \cite{\CSI, 4.1, 4.3})
is a direct sum of unipotent character sheaves on $G$; the set of all $A\in\hG^{un}$ which appear in 
$\ind_{P_J}^GA_0$ is denoted by $\hG^{un}_{J,\z}$ and is in canonical bijection with the set of isomorphism
classes of simple modules of $\End(\ind_{P_J}^GA_0)$ hence with $\Irr\WW^{\SS/J}$ (one can show that the
algebra $\End(\ind_{P_J}^GA_0)$ is canonically isomorphic to the group algebra of $\WW^{\SS/J}$).

\proclaim{Corollary 3.9} Assume that $G/\cz_G$ is simple or $\{1\}$. There exists a unique bijection 
$\fS_{\WW}@>\si>>\hG^{un}$ with the following property. If $A\in\hG^{un}$ corresponds to 
$(J,\e,\z)\in\fS_{\WW}$ then $A\in\hG^{un}_{J,\z}$ and $A$ corresponds as above to $\e\in\Irr\WW^{\SS/J}$.
\endproclaim

\subhead 3.10\endsubhead
We now return to the setup in 3.2. We assume that $G/\cz_G$ is simple or $\{1\}$. Combining the bijections 
in 3.5, 3.9 we obtain a canonical bijection $\hG^{un}\lra\fU_q$. Here are some properties of this bijection.
If $A\in\hG^{un}$ corresponds to $\r\in\fU_q$ then $\l_A\in\t_\r$;
moreover $(A:\ck_w)=(\r:R_w)$ for any $w\in W$.

\subhead 3.11 \endsubhead
Now assume that $\kk$ is as in 3.2 and that
$G$ has a fixed $\FF_p$-split rational structure with Frobenius map $F_0:G@>>>G$. We also assume that 
$F_0(B^*)=B^*$, $F_0(T^*)=T^*$. We fix a square root $\sqrt p$ of $p$ in $\bbq$. 
For any integer $t>0$, $F_0^t:G@>>>G$ is the Frobenius map for an $\FF_{p^t}$-rational structure on $G$
so that the definitions in 3.2 are applicable with $q=p^t$ and $F_0^t$ instead of $F$. (We set 
$\sqrt{q}=(\sqrt{p})^t$.)
For $w\in W$ we write $X_{w,t}$ instead of $X_w$ of 3.2 (thus $X_{w,t}=\{B\in\cb;(B,F_0^t(B))\in\fO_w\}$).
We want to show that in a certain sense the set $\hG^{un}$ is the limit of the sets $\fU_{p^t}$ as 
$t$ tends to $0$.

For $w\in\WW$ let $Y_{w,t}=\{(g,B)\in G\T\cb;(B,gF_0^t(B)g\i)\in\fO_w\}$ and let 
$\p_{w,t}:Y_{w,t}@>>>G$ be the first projection. Now $G$ acts on $Y_{w,t}$ by 
$x:(g,B)\m(xgF_0^t(x\i),xBx\i)$ and on $G$ (transitively) by $x:g\m xgF_0^t(x\i)$; then $\p_{w,t}$ is 
compatible with the $G$-actions. Let $\hG^{un}_t$ be the set of isomorphism classes of irreducible
perverse sheaves on $G$ which appear as
constituents of the perverse cohomology sheaf ${}^pH^i((\p_{w,t})_!\bbq)$ for some $w\in\WW,i\in\ZZ$. 
If $A\in\hG^{un}_t$ then $A$ is $G$-equivariant for a transitive $G$-action on $G$ hence it is of the
form $\cl[\dim G]$ where $\cl$ is an irreducible $G$-equivariant local system on $G$. Since the isotropy 
group at $1$ of the $G$-action on $G$ is $G^{F_0^t}$ we see that $\cl$ is completely determined by the 
(irreducible) representation of $G^{F_0^t}$ on the stalk at $1$ of $\cl$. This irreducible representation 
is clearly unipotent and we thus obtain a bijection $\hG^{un}_t\lra\fU_{p^t}$.
(Note that ${}^pH^i((\p_{w,t})_!\bbq)$ is up to shift the $G$-equivariant local system on $G$
such that the isotropy group $G^{F_0^t}$ acts on the stalk at $1$ as on the 
$G^{F_0^t}$-module $H^{i'}_c(X_{w,t},\bbq)$ for some $i'$.) On the other hand if we now take $t=0$ then 
$F_0^t$ becomes the identity map and $Y_{w,t},\p_{w,t}$ become $Y_w,\p_w$ in 2.1. 
The definition of $\hG^{un}_t$ specializes to the definition of $\hG^{un}$. 
In this sense the set $\hG^{un}$ can be viewed as limit of the sets $\fU_{p^t}$ as $t$ tends to $0$.

For $w\in W$ we define a map $\Ph_t:Y_{w,t}@>>>Y_{w,t}$ by $\Ph_t(g,B)=(g,gF_0^t(B)g\i)$.
This induces for any $i$ an isomorphism of perverse sheaves
${}^pH^i((\p_{w,t})_!\bbq)@>>>{}^pH^i((\p_{w,t})_!\bbq)$.
For any $\mu\in\bbq^*$ we denote by 
${}^pH^i((\p_{w,t})_!\bbq)_\mu$ the generalized $\mu$-eigenspace of the last isomorphism.
From the definitions we see that if $\t\in\fU_{p^t}$ corresponds as above to $A\in\hG^{un}_t$ then
we have $\mu\in\t_\r$ whenever $A$ is a constituent of 
${}^pH^i((\p_{w,t})_!\bbq)_\mu$ for some $w,i$.

When $t$ is replaced by $0$, the cosets of $\bbq^*$ modulo its subgroup $\{\sqrt{p}^{tr};r\in\ZZ\}$ become 
just the elements of $\bbq^*$ and for $w\in W$ the map
$\Ph_t:Y_{w,t}@>>>Y_{w,t}$ becomes the map $\Ph_0:Y_w@>>>Y_w$ given by $(g,B)=(g,gBg\i)$.
This induces for any $i$ an isomorphism of perverse sheaves
${}^pH^i(\p_{w!}\bbq)@>>>{}^pH^i(\p_{w!}\bbq)$ denoted again by $\Ph_t$.
For any $\mu\in\bbq^*$ we denote by 
${}^pH^i(\p_{w!}\bbq)_\mu$ the generalized $\mu$-eigenspace of the last isomorphism.
One can check from the definitions that if $A\in\hG^{un}$ is a constituent of
${}^pH^i(\p_{w!}\bbq)_\mu$ then $\mu=\l_A$. On the other hand $\l_A$ can be viewed as the limit as $t$ tends
to $0$ of the coset $\t_\r$ where $\r\in\fU_{p^t}$ corresponds to $A$.
		     
\subhead 3.12\endsubhead
Now assume that $G$ is the identity component of a possibly disconnected reductive group over $\kk$ and that
$G^1$ is a connected component of that reductive group. Then the notion of character sheaf on $G^1$ is
well defined (see \cite{\INT}.) Let $A$ be a character sheaf on $G^1$. Then to $A$ we can associate a root
of $1$ $\l_A\in\bbq^*$ as in 3.6. More precisely let $d$ be an integer $\ge1$ such that $g^d\in G$ for any 
$g\in G^1$. Note that $A$ is a simple perverse sheaf on $G^1$, $G$-equivariant for the conjugation action of
$G$ on $G^1$. Let $N$ be an open dense subset of the support of $A$ which is invariant by $G$-conjugation
and is such that $A|_N$ is up to shift a local system $\cl$. If $g\in N$ then the centralizer of $g$ in $G$ 
acts naturally and irreducibly 
on the stalk $\cl_g$. Since $g^d$ is in the centre of this centralizer it acts on $\cl_g$ as a nonzero
scalar which is independent of the choice of $N$ and is denoted by $\l_A$.

Using the invariant $\l_A$ we can find a combinatorial parametrization of the set of unipotent character
sheaves on $G^1$ parallel to that in 3.9.

Similarly, the set of unipotent representations of a not necessarily $\FF_q$-split connected group over 
$\FF_q$ can be parametrized in the same spirit using the invariant defined in \cite{\REP, 3.9}.

\enddocument